\documentclass[11pt,twoside]{amsart} 

\title{Equivariant Nielsen invariants for discrete groups}
\author{Julia Weber}

\address{Julia Weber\newline Max-Planck-Institut f\"ur Mathematik\newline Vivatsgasse 7\newline D-53111 Bonn\newline Germany\newline email: jweber@mpim-bonn.mpg.de}

\usepackage{bbm}
\usepackage{amssymb}
\usepackage{xspace}
\usepackage[all]{xy}\CompileMatrices
\usepackage{amsfonts}
\usepackage{amscd}

\SelectTips{cm}{}

\DeclareMathOperator{\IIm}{Im}
\DeclareMathOperator{\consub}{consub}
\DeclareMathOperator{\Mor}{Mor}
\DeclareMathOperator{\Aut}{Aut}
\DeclareMathOperator{\Ob}{Ob}
\DeclareMathOperator{\cst}{cst}

\DeclareMathOperator{\Deg}{Deg}
\DeclareMathOperator{\ch}{ch}
\DeclareMathOperator{\Fix}{Fix}
\DeclareMathOperator{\inc}{inc}
\DeclareMathOperator{\tr}{tr}
\DeclareMathOperator{\iind}{ind}
\DeclareMathOperator{\id}{id}

\DeclareMathOperator{\Is}{Is}
\DeclareMathOperator{\pr}{pr}

\newcommand{\Z}{\ensuremath{\mathbb{Z}}\xspace}
\newcommand{\Q}{\ensuremath{\mathbb{Q}}\xspace}
\newcommand{\Zz}{\ensuremath{\Z/2}\xspace}
\newcommand{\ind}[1]{\iind_{#1}}

\theoremstyle{plain}
\newtheorem{Def}{Definition}[section]
\newtheorem{Thm}[Def]{Theorem}

\newtheorem{Prop}[Def]{Proposition}
\newtheorem{Lem}[Def]{Lemma}
\newtheorem*{Thm*}{Theorem}

\theoremstyle{definition}

\newtheorem*{Def*}{Definition}
\newtheorem{Rem}[Def]{Remark}

\begin{document}

\begin{abstract}

For discrete groups $G$, we introduce equivariant Nielsen invariants. They are equivariant analogs of the
Nielsen number and give lower bounds for the number of fixed point
orbits in the $G$-homotopy class of an equivariant endomorphism $f\colon X\to X$. Under mild hypotheses,
these lower bounds are sharp. 

We use the equivariant Nielsen invariants to show that a
$G$-equi\-va\-riant endomorphism $f$ is $G$-homotopic to a fixed point free $G$-map if
the generalized equivariant Lefschetz invariant $\lambda_G(f)$ is zero. Finally, we prove a converse of the equivariant Lefschetz fixed
point theorem.\\



\end{abstract}
\maketitle

\section{Introduction}

The Lefschetz number is a classical invariant in algebraic topology. If the Lefschetz number $L(f)$ of an endomorphism $f\colon X\to X$ of a compact CW-complex is nonzero, then $f$ has a fixed point. This is stated by the famous Lefschetz fixed point theorem. The converse does not hold: If the Lefschetz number of $f$ is zero, then we cannot conclude $f$ to be fixed point free. 

A more refined invariant which allows to state the converse is the Nielsen number: The Nielsen number $N(f)$ is zero if and only if $f$ is homotopic to a fixed point free map. More generally, the Nielsen number is used to give precise minimal bounds for the number of fixed points of maps homotopic to $f$. Its development was started by Nielsen~\cite{nielsen}, a comprehensive treatment can be found in~\cite{jiang}.

We are interested in the equivariant generalization of these results. Given a discrete group $G$ and a $G$-equivariant endomorphism $f\colon X\to X$ of a finite proper $G$-CW-complex, we introduce equivariant Nielsen invariants
called $N_G(f)$ and $N^G(f)$. They are equivariant analogs of the
Nielsen number and are derived from the generalized equivariant Lefschetz invariant $\lambda_G(f)$~\cite[Definition~5.13]{weberunivlefschetz}. 

We proceed to show that these Nielsen invariants give minimal bounds for the number of
orbits of fixed points in the $G$-homotopy class of $f$. One even
obtains results concerning the type and ``location'' (connected
component of the relevant fixed point set) of these fixed point
orbits. These lower bounds are sharp if $X$ is a cocompact proper smooth $G$-manifold satisfying the standard gap hypotheses.

Finally, we prove a converse of the equivariant Lefschetz fixed
point theorem: If $X$ is a $G$-Jiang space as defined in Definition~\ref{defGJiangspace}, then $L_G(f)=0$ implies that $f$ is
$G$-homotopic to a fixed point free map. 
Here, $L_G(f)$ is the
equivariant Lefschetz class~\cite[Definition~3.6]{lueck-rosenberg},
the equivariant analog of the Lefschetz number.

These results were motivated by work of L\"uck and Rosenberg~\cite{lueck-rosenberg, lueck-rosenberg:equeulerchar}. For $G$ a discrete group and an endomorphism $f$ of a cocompact proper smooth $G$-manifold $M$, they prove an equivariant Lefschetz fixed point theorem~\cite[Theorem~0.2]{lueck-rosenberg}. The converse of that theorem is proven here.

Another motivation for the present article is the fact that the algebraic approach to the equivariant Reidemeister trace provides a good framework for computation. The connection to the machinery used in the study of transformation groups~\cite{lueck89, tomdieck} allows results to translate more readily from transformation groups to geometric equivariant topology and vice-versa. 

When $G$ is a compact
Lie group, Wong~\cite{wong93} obtains results on equivariant Nielsen numbers which strongly influenced us. The main difference between our
work and Wong~\cite{wong93} is that we treat possibly infinite discrete groups. Another difference is that our
approach is more structural. We can read off the
equivariant Nielsen invariants from the generalized equivariant
Lefschetz
 invariant $\lambda_G(f)$.

In case $G$ is a finite group, Ferrario~\cite{ferrario03moebius} studies a collection of generalized Lefschetz numbers which can be thought of as an equivariant generalized Lefschetz number. In contrast to the generalized equivariant
Lefschetz invariant $\lambda_G(f)$ these do not incorporate the $WH$-action on the fixed point set $X^H$. A generalized Lefschetz trace for equivariant maps has also been defined by Wong~\cite{wongprivate}, as mentioned in~\cite{hart}. 

For compact Lie groups, earlier definitions of generalized Lefschetz numbers for equivariant maps were made by~\cite{wilczinski, fadell-wong}. They used the collection of generalized Lefschetz numbers of the maps $f^H$, for $H<G$. In general, these numbers are not sufficient since they do not take the equivariance into account adequately. For further reading on equivariant fixed point theory, see~\cite{ferrario05}, where an extensive list of references is given.\\

This paper is organized as follows. In Section~\ref{sec2}, we introduce the generalized equivariant Lefschetz invariant. We briefly assemble the concepts and definitions which are needed for the definition of equivariant Nielsen invariants in Section~\ref{sec3}. 

The equivariant Nielsen invariants give lower
bounds for the number of fixed point orbits of
$f$. This is shown in Section~\ref{sec4}. The standard gap hypotheses are introduced, and it is shown that under these hypotheses these lower bounds are sharp.

In the non-equivariant case, we know that the generalized Lefschetz invariant is the right
element when looking for a precise count of fixed points. We read off
the Nielsen number from this invariant. In general, the Lefschetz
number contains too little information. 

But under certain conditions, we can
conclude facts about the Nielsen number from the Lefschetz number
directly. These conditions are called Jiang conditions~\cite[Definition~II.4.1]{jiang}~\cite[Chapter~VII]{brown}. In Section~\ref{sec5}, we introduce the equivariant version of these conditions. We also give examples of $G$-Jiang spaces.

In Section~\ref{sec6}, we derive equivariant analogs of statements about Nielsen
numbers found in Jiang~\cite{jiang}, generalizing results of Wong~\cite{wong93} to infinite
discrete groups. In particular, if $X$ is a $G$-Jiang space, the
converse of the equivariant Lefschetz fixed point theorem holds.

\section{The generalized equivariant Lefschetz invariant}\label{sec2}

Classically, the Nielsen number is defined geometrically by counting essential fixed
point classes~\cite[Chapter~VI]{brown}~\cite[Definition~I.4.1]{jiang}. Alternatively one defines it using the generalized
Lefschetz invariant.

Let $X$ be a finite CW-complex, let $f\colon X\to
X$ be an endomorphism, let $x$ be a basepoint of $X$, and
let 
\[
\lambda(f)=\sum_{\overline{\alpha}\in \pi_1(X,x)_{\phi}}
n_{\overline{\alpha}} \cdot \overline{\alpha} \in \Z \pi_1(X,x)_{\phi} 
\]
be the \emph{generalized Lefschetz invariant} associated to $f$~\cite{reidemeister,wecken2}, where 
\[
\Z \pi_1(X,x)_{\phi}:=\Z \pi_1(X,x)/\phi(\gamma)\alpha\gamma^{-1}\sim \alpha, \text{ with } \gamma,\alpha\in \pi_1(X,x).
\]
Here $\phi$ is the map induced by $f$ on the fundamental group $\pi_1(X,x)$. We have $\phi(\gamma)=wf(\gamma)w^{-1}$, where $w$ is a path from $x$ to $f(x)$. The generalized Lefschetz invariant is also called \emph{Reidemeister trace} in the literature, which goes back to the original name "Reidemeistersche Spureninvariante" used by Wecken~\cite{wecken2} for the invariant constructed by Reidemeister~\cite{reidemeister}.

The set $\pi_1(X,x)_{\phi}$ is often denoted by $\mathcal{R}(f)$ and called set of \emph{Reidemeister classes} of $f$. Since we will introduce a variation of this set in Definition~\ref{defzgroup}, we prefer to stick with the notation used in~\cite{weberunivlefschetz}.

\begin{Def}\label{defnielsennumber}
The \emph{Nielsen number} of $f$ is defined by
\[
N(f):=\# \bigl\{\overline{\alpha}\; \big| \;  n_{\overline{\alpha}}\neq 0 \bigr\}.
\]
\end{Def}

The Nielsen number is the number of classes in
$\pi_1(X,x)_\phi$ with non-zero coefficients. A class $\alpha$ with
non-zero coefficient corresponds to an essential fixed point class in
the geometric sense.

In the equivariant setting, the fundamental category replaces the fundamental group. The fundamental category of a topological
space $X$ with an action of a discrete group $G$ is defined as follows~\cite[Definition~8.15]{lueck89}. 

\begin{Def}\label{deffundcat}
Let $G$ be a discrete group, and let $X$ be a $G$-space. Then the
\emph{fundamental category} $\Pi(G,X)$ is the following category:
\begin{itemize}
\item The objects $\Ob\bigl(\Pi(G,X)\bigr)$ are $G$-maps $x\colon G/H\to X$,
      where the $H\leq G$ are subgroups. 
\item The morphisms $\Mor\bigl(x(H),y(K)\bigr)$ are pairs $(\sigma,[w])$, where
 \begin{itemize}
   \item $\sigma$ is a $G$-map $\sigma\colon G/H\to G/K$
     \item $[w]$ is a homotopy class of $G$-maps $w\colon G/H\times I
           \to X$ relative $G/H\times \partial I$ such that $w_1=x$
           and $w_0=y\circ \sigma$. 
\end{itemize}
\end{itemize}
\end{Def}

The fundamental category is a combination of the orbit category of $G$
and the fundamental groupoid of $X$. If $X$ is a point, then the
fundamental category is
just the orbit category of $G$, whereas when $G$ is the trivial group,
the definition reduces to the definition of the fundamental groupoid
of $X$. 

We often view $x$ as the point $x(1H)$ in the fixed point set $X^H$. We
call $X^H(x)$ the connected component of $X^H$ containing $x(1H)$. We also consider the relative fixed point set, the pair $\bigl(X^H(x),X^{>H}(x)\bigr)$. Here $X^{>H}(x)=\{z\in X^H(x) \; | \; G_z\neq H\}$ is the singular set, where $G_z$ denotes the isotropy group of $z$. In order to simplify notation, we use $f^H(x)$ to denote
$f|_{X^H(x)}$, and we use $f_H(x)$ instead of
$f|_{\bigl(X^H(x),X^{>H}(x)\bigr)}$. 

Fixed points of $f$ can only exist in $X^H(x)$ when $X^H(f(x))=X^H(x)$, i.e, when the points $f(x)$ and $x$ lie in the same connected component of $X^H$.  

\begin{Def}\label{defzgroup}
For $x\in \Ob\Pi(G,X)$ with $X^H(f(x))=X^H(x)$ and a morphism
$v=(\id,[w])\in \Mor(f(x),x)$, set
\[
\Z \pi_1\bigl(X^H(x),x\bigr)_{\phi'}:=
\Z\pi_1\bigl(X^H(x),x\bigr)/\phi(\gamma)\alpha \gamma^{-1}\sim \alpha,
\]
where $\alpha\in \pi_1(X^H(x),x)$, $\gamma\in \Aut(x)$ and
$\phi(\gamma)=v\phi(\gamma)v^{-1}\in \Aut(x)$.
\end{Def}

The automorphism group $\Aut(x)$ of the object $x$ in the category $\Pi(G,X)$ is a group extension lying in the short exact sequence
\[
1\to \pi_1(X^H(x),x)\to \Aut(x)\to WH_x \to 1.
\]
Here $WH:=N_G H/H$ is the Weyl group of $H$, it acts on $X^H$. We call $WH_x$ the subgroup of $WH$ which fixes the connected component $X^H(x)$. 

Groups obtained from different choices of the path $w$ and of the point $x$ in its isomorphism class $\overline{x}$ are canonically isomorphic, so those choices do not play a role. The group $\Z \pi_1\bigl(X^H(x),x\bigr)_{\phi'}$ generalizes the group $\Z \pi_1(X,x)_{\phi}$ defined above. So it can be seen as the free abelian group generated by equivariant Reidemeister classes of $f^H(x)$ with respect to the action of $WH_x$ on $X^H(x)$. 

A map $(\sigma,[w])\in \Mor(x,y)$ induces a group homomorphism 
\[
(\sigma,[w])^*\colon \Z \pi_1\bigl(X^K(y),y\bigr)_{\phi'}\to \Z \pi_1\bigl(X^H(x),x\bigr)_{\phi'}
\]
by twisted conjugation, and we know that the induced group homomorphism is the same for every map in $\Mor(x,y)$~\cite[Lemma~5.2]{weberunivlefschetz}.

The \emph{generalized equivariant Lefschetz invariant}~\cite[Definition~5.13]{weberunivlefschetz}, $\lambda_G(f)$, is an element in the group
\[
\Lambda_G(X,f):=\bigoplus_{\overline{x}\in \Is \Pi(G,X),
  \atop X^H(f(x))=X^H(x)} \Z \pi_1(X^H(x),x)_{\phi'}.
\]
Here $\Is \Pi(G,X)$ denotes the set of isomorphism classes of the category $\Pi(G,X)$. Geometrically, it corresponds to the set of $WH$-orbits of connected components $X^H(x)$ of the fixed point sets $X^H$, for $(H)\in \consub(G)$, i.e., for a set of representatives of conjugacy classes of subgroups of $G$. There is a bijection $\Is \Pi(G,X)\xrightarrow{\simeq} \amalg_{(H)\in \consub(G)} WH\setminus \pi_0(X^H)$ which sends $x\colon G/H\to X$ to the orbit under the WH-action on $\pi_0(X^H)$ of the component $X^H(x)$ of $X^H$ which contains the point $x(1H)$~\cite[Equation~3.3]{lueck-rosenberg}.

Let $\widetilde{f^H(x)}$ and $\widetilde{f^{>H}(x)}$ denote the lift of $f^H(x)$ to the universal covering space $\widetilde{X^H(x)}$ and to the subset $\widetilde{X^{>H}(x)}\subseteq\widetilde{X^H(x)}$ that projects to $X^{>H}(x)$ under the covering map.

At the summand indexed by $\overline{x}$, the generalized equivariant Lefschetz invariant is given by 
\[
\lambda_G(f)_{\overline{x}} :=L^{\Z \Aut(x)}\bigl(\widetilde{f^H(x)},\widetilde{f^{>H}(x)}\bigr)\in \Z\pi_1(X^H(x),x)_{\phi'},
\] 
where the \emph{refined equivariant Lefschetz number}~\cite[Definition~5.7]{weberunivlefschetz} appears on the right hand side. It is defined by
\[
L^{\Z \Aut(x)}\bigl(\widetilde{f^H(x)},\widetilde{f^{>H}(x)}\bigr)  :=  \sum_{p\geq 0}(-1)^p \tr_{\Z \Aut(x)}(C^c_p(\widetilde{f^H(x)},\widetilde{f^{>H}(x)})),
\]
where the trace map $\tr_{\Z \Aut(x)}$~\cite[Definition~5.4]{weberunivlefschetz} 
is induced by the projection $\Z \Aut(x)\to \Z\pi_1(X^H(x),x)_{\phi'}, \sum_{g\in \Aut(x)}r_g\cdot g\mapsto \sum_{g\in \pi_1(X^H(x),x)} r_g\cdot \overline{g}$. Instead of $\Z$, other rings can be used. 

This trace map generalizes the trace map used in~\cite{lueck-rosenberg}, and the refined equivariant Lefschetz number is a generalization of the orbifold Lefschetz number~\cite[Definition~1.4]{lueck-rosenberg}.

The refined equivariant Lefschetz number $L^{\Q \Aut(x)}\bigl(\widetilde{f^H(x)}\bigr)$ will be particularly important to us, so we give some formulas describing it. For a finite proper $G$-CW-complex $X$ we have~\cite[Lemma~5.9]{weberunivlefschetz} 
\[
L^{\Q \Aut(x)}\bigl(\widetilde{f^H(x)}\bigr)
= \sum_{p\geq 0} (-1)^p \hspace{-1em} \sum_{G\cdot e\in G \setminus I_p(X)} \hspace{-1em} |G_e|^{-1} \cdot \inc_{\phi}(f,e) \in \Q \pi_1(X^H(x),x)_{\phi'}.
\]
Here $I_p(X)$ denotes the set of $p$-cells of $X$, $e$ runs through the equivariant cells of $X$, and $G_e$ is its isotropy group. The \emph{refined incidence number}~\cite[Definition~5.8]{weberunivlefschetz} $\inc_\phi(f,e)\in \Z\pi_1(X^H(x),x)_{\phi'}$ for a $p$-cell $e\in I_p(X)$ is defined to be the "degree" of the composition
\begin{align*}
\overline{e}/\partial e & \xrightarrow{i_e} \hspace{-1em}\bigvee_{e'\in I_p(X)}
\overline{e'}/\partial e' \xrightarrow{h \sim} X_p/X_{p-1}
\xrightarrow{f} X_p/X_{p-1} \xrightarrow{h^{-1} \sim} \hspace{-1em} \bigvee_{e'\in
  I_p(X)} \overline{e'}/\partial e'\\
& \xrightarrow{\pr_{\pi\cdot \overline{e}/\partial e }} \pi \cdot
\overline{e}/\partial e \xrightarrow{\overline{\, \cdot \,}}\pi_{\phi'} \cdot
\overline{e}/\partial e.
\end{align*}
Here $\overline{e}$ is the closure of the open $p$-cell $e$ and $\partial e=\overline{e}\setminus e$. The map $i_e$ is the inclusion, $h$ is a homeomorphism and $\pr_{\pi\cdot \overline{e}/\partial e }$ is the projection.

If $X=M$ is a cocompact proper $G$-manifold, we have~\cite[Theorem~6.6]{weberunivlefschetz}  
\[
L^{\Q \Aut(x)}\bigl(\widetilde{f^H(x)}\bigr)
=  \hspace{-1.5em} \sum_{WH_x \cdot z\in \atop WH_x \setminus
  \Fix (f^H(x) )} \hspace{-1.5em} \bigl|(WH_x)_z\bigr|^{-1} \deg\bigl(
\bigl(\id_{T_z M^H(x)} - T_z(f^H(x))\bigr)^c  \bigr) \cdot \overline{\alpha_z}.
\]
Here the map on the tangent space is extended to the one-point compactification $\bigl(T_z M^H(x) \bigr)^c$. The relative versions of these formulas also hold.

We have $L^{\Q
  \Aut(x)}\bigl(\widetilde{f^H(x)}\bigr)=\ch_G(X,f)\bigl(\lambda_G(f)\bigr)_{\overline{x}}$~\cite[Lemma~6.4]{weberunivlefschetz}, where $\ch_G(X,f)\colon \Lambda_G(X,f)\to \bigoplus_{\overline{y}\in \Is\Pi(G,X)} \Q \pi_1(X^K(y),y)_{\phi'}$ is the character map~\cite[Definition~6.2]{weberunivlefschetz}. So we can derive $L^{\Q
  \Aut(x)}\bigl(\widetilde{f^H(x)}\bigr)$ from $\lambda_G(f)$.

The equivariant analog of the Lefschetz number is the \emph{equivariant Lefschetz class} $L_G(f)\in \bigoplus_{\overline{x}\in \Is \Pi(G,X),
  X^H(f(x))=X^H(x)} \Z $, whose value at $\overline{x}$ is given by $L_G(f)_{\overline{x}}=L^{\Z WH_x}\bigl({f^H(x)},{f^{>H}(x)}\bigr)$~\cite[Definition~3.6]{lueck-rosenberg}. The projection of $\pi_1(X^H(x),x)_{\phi'}$ to the trivial group $\{1\}$ induces an augmentation map 
\[
s\colon \hspace{-0.5em}\bigoplus_{\overline{x}\in \Is \Pi(G,X),
  \atop X^H(f(x))=X^H(x)} \hspace{-1.5em} \Z \pi_1(X^H(x),x)_{\phi'}\hspace{0.5em}\to\hspace{-0.5em}\bigoplus_{\overline{x}\in \Is \Pi(G,X),
  \atop X^H(f(x))=X^H(x)} \hspace{-1.5em} \Z 
\]
which sends $\lambda_G(f)$ to $L_G(f)$.

\section{Equivariant Nielsen Invariants}\label{sec3}

In this section, we define equivariant Nielsen invariants. Given an
element $\sum_{\overline{\alpha}} n_{\overline{\alpha}}\cdot \overline{\alpha} \in \Z
\pi_1(X^H(x),x)_{\phi'}$, we call a class $\overline{\alpha}\in
\pi_1(X^H(x),x)_{\phi'}$ \emph{essential} if the coefficient
$n_{\overline{\alpha}}$ is non-zero. 

Let $G$ be a discrete group and let $X$ be a cocompact proper smooth
$G$-manifold. Let $f\colon X\to X$ be a smooth $G$-equivariant
map such that $\Fix(f)\cap \partial X = \emptyset$ and such that for
every $z\in \Fix(f)$ the determinant of the map $(\id_{T_z X} - T_z
f)$ is different from zero. One can always find a representative in the $G$-homotopy class of $f$ which satisfies this assumption. Since the generalized equivariant Lefschetz invariant is $G$-homotopy invariant, we can replace $f$ by this representative if necessary. 

\begin{Def}
The \emph{equivariant Nielsen class} of $f$ is
\begin{eqnarray*}
\nu_G(f) & = & \sum_{Gz\in G\setminus \Fix(f)} \frac{\det\bigl(\id_{T_z X}-T_z(f)\bigr)}{\left|\det\bigl(\id_{T_z X}-T_z(f)\bigr)\right|} \cdot \overline{\alpha_z} .
\end{eqnarray*}
\end{Def}

Here ${\alpha_z}\in \pi_1(X^{G_z}(z),x)$ is the loop given by $[t*f(t)^{-1}*w]$, where $x$ is a basepoint in $X^{G_z}(z)$, $t$ is a path from $x$ to $z$ and $w$ is a path from $f(x)$ to $x$. The basepoint $x$ may differ from $z$, e.g., if we have more than one fixed point in a connected component of $X^{G_z}$. If $x=z$, we may choose $t$ and $w$ to be constant. The equivalence relation assures that this definition is independent of the choices involved.

We can also derive the equivariant Nielsen class $\nu_G(f)$ from the generalized equivariant Lefschetz invariant $\lambda_G(f)$. 

\begin{Lem} \label{lemnulambda}
The invariant $\nu(f)$ is the image of the generalized equivariant Lefschetz
invariant $\lambda(f)$ under the quotient map where we divide out the images of
non-isomorphisms. 
\[
\nu_G(f)=\overline{\lambda_G(f)}\in \bigoplus_{\overline{x}\in \Is \Pi(G,X),
  \atop X^H(f(x))=X^H(x)} \hspace{-2em} \Z \pi_1(X^H(x),x)_{\phi'}/\{\IIm (\sigma,[w])^*\; | \; \sigma \text{ non-isom.}\}
\]
\end{Lem}

\begin{proof}
We consider the equation obtained in the
refined equivariant Lefschetz fixed point theorem~\cite[Theorem~0.2]{weberunivlefschetz}. We have
\[
\lambda_G(f)=\sum_{G z \in G \setminus \Fix(f)}
\Lambda_G(z,f)\circ \ind{G_z\subseteq G} (\Deg_0^{G_z}((\id_{T_z X} -
  T_z f)^c)).
\]
Here $\Deg_0^{G_z}$ is the equivariant degree~\cite{lueck-rosenberg}, it has values in the Burnside ring $A(G_z)$. On basis elements $[G_z/L]\in A(G_z)$, the map $\Lambda_G(z,f)\circ \ind{G_z\subseteq G}$ is given by 
\[
\Lambda_G(z,f)\circ \ind{G_z\subseteq G}([G_z/L])= (\pr,[\cst])^* \overline{\alpha_z},
\]
where $\cst$ denotes the constant map and $(\pr,[\cst])^*\colon \Z \pi_1(X^{G_z}(z),x)_{\phi'}\to \Z \pi_1(X^L(z\circ \pr),x\circ \pr)_{\phi'}$ is the map induced by the projection $\pr\colon G_z/L \to G_z/G_z$.

We know that $\Deg_0^{G_z}((\id_{T_z X} -
  T_z f)^c)$ is a unit of the Burnside ring $A(G_z)$ since $\bigl(\Deg_0^{G_z}((\id_{T_z X} -
  T_z f)^c)\bigr)^2=1$~\cite[Example~4.7]{lueck-rosenberg}. In general a unit of the Burnside ring $A(G_z)$ may consist of
  more than one summand~\cite{tomdieck79}. The summand $[G_z/G_z]$ is
  always included with a coefficient $+1$ or $-1$, but there might be
  summands $[G_z/L]$ for $L<G_z$ appearing. So one fixed point might give
  more than one class with non-zero coefficients. 

If we divide out the images of non-isomorphisms, then we divide out the image of $(\pr,[\cst])^*$ for all $L\neq G_z$. We are left with the summand~$\pm \overline{\alpha_z}$ coming from $\pm 1
  [G_z/G_z]$. This cannot lie in the image of any non-isomorphism. So each fixed
  point leads to exactly one summand. The sign is the sign of the determinant $\det\bigl(\id_{T_z X}-T_z(f)\bigr)$, so the claim follows.
\end{proof}

We set 
\[
\Z \pi_1(X^H(x),x)_{\phi''}:=\Z \pi_1(X^H(x),x)_{\phi'}/\{\IIm (\sigma,[w])^*\; | \; \sigma \text{ non-isom.}\}.
\]
We use the equation established in Lemma~\ref{lemnulambda} to define $\nu_G(f)$ directly for all endomorphisms of finite proper $G$-CW-complexes.

\begin{Def}
Let $X$ be a finite proper $G$-CW-complex, and let $f\colon X\to X$ be an equivariant endomorphism.
Then the \emph{equivariant Nielsen class} of $f$ is
\[
\nu_G(f):=\overline{\lambda_G(f)}\in \bigoplus_{\overline{x}\in \Is \Pi(G,X),
  \atop X^H(f(x))=X^H(x)} \hspace{-1em} \Z \pi_1(X^H(x),x)_{\phi''}.
\]
\end{Def}

We define equivariant Nielsen invariants by counting the essential
 classes $\overline{\alpha}$ of ${\nu_G(f)}_{\overline{x}}$ in $\Z
 \pi_1(X^H(x),x)_{\phi''}$ and
 of $L^{\Q \Aut(x)}\bigl(\widetilde{f^H(x)}\bigr)$ in $\Q \pi_1(X^H(x),x)_{\phi'}$. 

\begin{Def}
Let $G$ be a discrete group, let $X$ be a finite proper $G$-CW-complex,
and let $f\colon X\to X$ be a $G$-equivariant map.  

Then the \emph{equivariant Nielsen
  invariants} of $f$ are elements 
\[
N_G(f), N^G(f) \in \bigoplus_{\overline{x}\in \Is \Pi(G,X)} \Z
\]
defined for $\overline{x}$ with $X^H(f(x))=X^H(x)$ by
\begin{eqnarray*}
{N_G(f)}_{\overline{x}} & := & \hspace{-1ex} \# \bigl\{
\text{essential classes of }\nu_G(f)_{\overline{x}}\bigr\}\\
{N^G(f)}_{\overline{x}} & := &  \hspace{-1ex} \min \Bigl\{ \# \mathcal{C}\; \Big| \;
\mathcal{C}\subseteq \bigcup_{y\geq x} \pi_1(X^K(y),y)_{\phi'} \text{
  such that for all } \overline{z} \geq \overline{x}\text{ and} \\
& & \hspace{-2ex} \qquad \text{for all essential classes $\overline{\alpha}$ of }
  L^{\Q \Aut(z)}\bigl(\widetilde{f^{G_z}(z)}\bigr)\text{ there are}
  \\
& & \hspace{-2ex} \qquad\overline{\beta}\in
  \mathcal{C} \text{ and } (\sigma,[t])\in \Mor(z,y_{\overline{\beta}})\text{ such that~}
  (\sigma,[t])^*(\overline{\beta})=\alpha
  \Bigr\} .
\end{eqnarray*}
We continue them by $0$ to $\overline{x}\in \Is \Pi(G,X)$ with $X^H(f(x))\neq X^H(x)$.
\end{Def}

Note that
${N_G(f)}_{\overline{x}}=N_G(f_H(x))$ and
${N^G(f)}_{\overline{x}}=NO_G(f^H(x))$ in the notation of Wong~\cite{wong93}. Thus the invariants defined here using the algebraic approach are equivalent to the invariants defined using the classical covering space approach of Wong.

An essential class $\overline{\alpha}$  of
$L^{\Q
  \Aut(x)}\bigl(\widetilde{f^H(x)}\bigr)$ corresponds to an essential fixed
point class of $f^H(x)$, a $WH_x$-orbit of fixed points which one
cannot get rid of under any $G$-homotopy, as can be seen from the refined orbifold
Lefschetz fixed point theorem~\cite[Theorem~6.6]{weberunivlefschetz}. An essential class
$\overline{\alpha}$ of $\nu_G(f)_{\overline{x}}$ corresponds to an
essential fixed point class of $f_H(x)$, an orbit of fixed points on
$X^H(x)\setminus X^{>H}(x)$ that cannot be moved into
$X^{>H}(x)$. Counting the essential classes will give us
information on the number of fixed points and fixed point orbits.  

The equivariant Nielsen invariants are $G$-homotopy invariant since they are derived
from $\lambda_G(f)$, which is itself $G$-homotopy invariant.

\begin{Prop}
Given a $G$-homotopy $f\simeq_G f'$, we have 
\begin{eqnarray*}
N_G(f) & = & N_G(f')\\
N^G(f) & = & N^G(f').
\end{eqnarray*}
\end{Prop}

\begin{proof}
If $f\simeq_G f'$, with a homotopy $H\colon X\times I\to X$ such that
$H_0=f$ and $H_1=f'$, then by invariance under homotopy equivalence~\cite[Theorem~5.14]{weberunivlefschetz} we have an isomorphism $\Lambda_G(i_1)^{-1}
\Lambda_G(i_0)\colon \Lambda_G(X,f)\xrightarrow{\sim} \Lambda_G(X,f')$
which sends $\lambda_G(f)$ to $\lambda_G(f')$. The isomorphisms
$\Lambda_G(i_1)$ and $\Lambda_G(i_0)$ are given by composition of
maps, so they do not change the number of essential classes. They also do not change the property of a class to lie in
the image of a non-isomorphism. So we have $N_G(f)=N_G(f')$.

An isomorphism ${i_0}_*
\colon \Q \Pi(G,X)_{\phi,\overline{y}}\to \Q \Pi(G,X\times
I)_{\Phi,\overline{i_0(y)}}$ is induced by the inclusion $i_0$, and
analogously $i_1$ induces an isomorphism. These
isomorphisms do not change the number of essential classes. We have
$\ch_G(X,f)(\lambda_G(f))=({i_0}_*)^{-1}{i_1}_*\ch_G(X,f')(\lambda_G(f'))$,
so $N^G(f)=N^G(f')$. 
\end{proof}

\section{Lower Bound Property}\label{sec4}

The equivariant Nielsen invariants give a lower
bound for the number of fixed point orbits on $X^H(x)\setminus X^{>H}(x)$ and on $X^H(x)$, for maps lying in the $G$-homotopy class of
$f$. Under mild hypotheses, this is even a sharp lower bound.  

\begin{Def}
Let $G$ be a discrete group, let $X$ be a finite proper $G$-CW-complex, and let $f\colon X\to X$
be a $G$-equivariant map. For every $\overline{x}\in\Is\Pi(G,X)$, with
$x\colon G/H\to X$, we set
\begin{eqnarray*}
{M_G(f)}_{\overline{x}} & := & \min \bigl\{\#
\text{ fixed point orbits of }\varphi_H(x) \, \big| \,\varphi\simeq_G f \bigr\},\\
M^G(f)_{\overline{x}} & := & \min \bigl\{ \# \text{ fixed point orbits of
}\varphi^H(x) \, \big| \,\varphi\simeq_G f\bigr\}.
\end{eqnarray*}
\end{Def}

When speaking of fixed point orbits of $f^H(x)$, we can either look at the $WH_x$-orbits $WH_x\cdot
z\subseteq X^H(x)$ or at the $G$-orbits $G \cdot z\subseteq
X^{(H)}(x)$, for a fixed point $z$ in $X^H(x)$. These two notions are
of course equivalent. 

We now proceed to show the first important property of the equivariant
Nielsen invariants, the lower bound property.

\begin{Prop} \label{proplowerbound}
For every $\overline{x}\in \Is\Pi(G,X)$ we have
\begin{eqnarray*}
{N_G(f)}_{\overline{x}} & \leq & {M_G(f)}_{\overline{x}}\\
N^G(f)_{\overline{x}} & \leq & M^G(f)_{\overline{x}}.
\end{eqnarray*}
\end{Prop}

\begin{proof}
1)
If $\overline{\alpha}\in \pi_1(X^H(x),x)_{\phi''}$ is an essential
class of $\nu_G(f)_{\overline{x}}$, then there has to be at least one
fixed point orbit in $X^H(x)\setminus X^{>H}(x)$ that corresponds to $\overline{\alpha}$ and that
cannot be moved into $X^{>H}(x)$. So, for any $\varphi\simeq_G f$, the
restriction $\varphi_{H}$ must have
at least ${N_G(f)}_{\overline{x}}$ fixed point orbits in
$X^{H}(x)\setminus X^{>H}(x)$. We arrive at ${N_G(f)}_{\overline{x}}
\leq \{ \# \text{fixed point
  orbits of }\varphi_H \}$ for all
$\varphi\simeq_G f$, so ${N_G(f)}_{\overline{x}}\leq {M_G(f)}_{\overline{x}}$.

2)
Let $\overline{x}\in \Is \Pi(G,X)$. Suppose that $\varphi\simeq_G
f$ such that $\varphi^H(x)$ has $M^G(f)_{\overline{x}}$ fixed point orbits in
$X^H(x)$. Let $\mathcal{C}\subseteq \bigcup_{\overline{x}\leq \overline{y}}
\pi_1(X^K(y),y)_{\phi'}$ such that
$N^G(\varphi)_{\overline{x}}=N^G(f)_{\overline{x}}=\#\mathcal{C}$. If there
were less than $\# \mathcal{C}$ fixed point orbits in $X^H(x)$, there would be
less that $\# \mathcal{C}$ essential classes and we could have
chosen a smaller $\mathcal{C}$. So there are at least $\#\mathcal{C}$
essential classes, and thus $\varphi^H(x)$
has at least $\#\mathcal{C}$ fixed point orbits.  
\end{proof}

To prove the sharpness of this lower bound, we need certain
hypotheses, which are usually introduced when dealing with these
problems. Such conditions were first used in~\cite{fadell-wong}. Some authors treat slightly weakened
assumptions~\cite{ferrario99, ferrario03moebius, jezierski,
  wilczinski}. We do not weaken the standard gap hypotheses in the
context of functorial equivariant Lefschetz
invariants since the standard gap hypotheses are not homotopy
invariant. So an analog of Theorem~\ref{thmconvequL} would not hold.

\begin{Def}\label{defstandardhyp}
Let $G$ be a discrete group and let $X$ be a cocompact smooth
$G$-manifold. We say that $X$ satisfies the \emph{standard gap
  hypotheses} if for each
$\overline{x}\in \Is \Pi(G,X)$, with $x\colon G/H \to X$, the inequalities $\dim X^H(x)\geq 3$
and $\dim X^H(x) - \dim X^{>H}(x) \geq 2$ hold.  
\end{Def}

Under these hypotheses, we can use an equivariant analog of the
classical Wecken method~\cite{wecken2} to coalesce fixed points.  

\begin{Lem} \label{lemcoalesce}
Let $G$ be a discrete group and let $X$ be a cocompact proper smooth $G$-manifold
satisfying the standard gap hypotheses. Let $f\colon X\to X$ be
a $G$-equivariant map. Let $\mathcal{O}_1=G x_1$ and
$\mathcal{O}_2 = G x_2$ be two distinct isolated $G$-fixed
point orbits, where $x_1:G/H\to X$ and $x_2:G/K \to X$
with $x_1\leq x_2$. 
Suppose that there are paths
$(\sigma_1,[t_1])\in \Mor(x,x_1)$ and $(\sigma_2,[t_2])\in
\Mor(x,x_2)$ for an $\overline{x}\in \Pi(G,X)$, with $x\colon G/H \to X$, such that
$(\sigma_1,[t_1])^*\overline{1_{x_1}}=\overline{\alpha}=(\sigma_2,[t_2])^*\overline{1_{x_2}}$,
i.e., that the fixed point orbits induce the same
$\overline{\alpha}\in \pi_1(X^H(x),x)_{\phi'}$.
 Then there exists a $G$-homotopy
$\{f_t\}$ relative to $X^{>(H)}$ such that $f_0=f$ and $\Fix f_1=\Fix
f_0 - G \mathcal{O}_1$. 
\end{Lem}

\begin{proof}
Suppose first that $\overline{x_1}< \overline{x_2}$. Then
$\Mor(x_1,x_2)\neq \emptyset$. By replacing
$x_1$ and $x_2$ with other points in the orbit if necessary, we can
suppose that there exists a morphism
$(\tau,[v])\in\Mor(x_1,x_2)$, where $v$ is a path in
$X^H(x)$  with $v_1=x_1$
and $v_0=x_2\circ \tau$ and $\tau\colon G/H\to G/K$ is a projection.
We know that $v\simeq f^H\circ v$ (relative endpoints). (This is an
equivalent characterization of $x_1$ and $x_2$ belonging to the same
fixed point class~\cite[I.1.10]{jiang}.) 
Since
$x_1\in X^H(x)\setminus X^{>H}(x)$ and $x_2\in X^{>H}(x)$ and $\dim X^H(x) -
\dim X^{>H}(x) \geq
2$, we may assume that $v$ can be chosen such that $v((0,1])\subseteq
X^H(x)\setminus X^{>H}(x)$. We coalesce $x_1$ and $x_2$ along $v$ as
in~\cite[1.1]{wong91_location},~\cite[6.1]{schirmer}. 
We can do this by only changing $f$ in a
(cone-shaped) neighborhood $U(v)$ of $v$. Because of the proper action of $G$
on $X$ and the free action of $WH$ on $X^H\setminus X^{>H}$, this neighborhood
$U(v)$ can be chosen such that in $X^H\setminus X^{>H}$ it does not
intersect its $g$-translates for $g\not\in H\leq G$. Taking the
$G$-translates of $U(v)$, we move $\mathcal{O}_1$ to $\mathcal{O}_2$
along the paths $Gv$ in $GU(v)$, not changing the map $f$ outside $GU(v)$.  

Now suppose $\overline{x_1}=\overline{x_2}$. In this case, the result
follows from the result of Wong~\cite[5.4]{wong91_equivnielsen} since
$X^H(x)\setminus X^{>H}(x)$ is a free and proper $WH_x$-space, where
again the proper action of $G$ on $X$ ensures that we can
find a neighborhood of a path from $x_1$ to $x_2$ such that the
$G/H$-translates do not intersect. 
\end{proof}

From Lemma~\ref{lemcoalesce}, we can conclude the sharpness of the lower bound
given by the equivariant Nielsen invariants. 

\begin{Thm} \label{thmdoubtful}
Let $G$ be a discrete group. Let $X$ be a cocompact proper smooth
$G$-manifold satisfying the standard gap hypotheses. Let
$f\colon X\to X$ be a $G$-equivariant endomorphism. Then 
\begin{eqnarray*}
M_G(f)_{\overline{x}} & = & N_G(f)_{\overline{x}}\\
M^G(f)_{\overline{x}} & = & N^G(f)_{\overline{x}}
\end{eqnarray*}
for all $\overline{x}\in \Is \Pi(G,X)$. 
\end{Thm}

\begin{proof}
1)
Since $X$ is a cocompact smooth $G$-manifold, there is a $G$-map $f'$
which is $G$-homotopic to $f$ and which has only finitely
many fixed point orbits. We apply Lemma~\ref{lemcoalesce} to $f'$ to
coalesce fixed point orbits in $X^H(x)\setminus X^{>H}(x)$ with others of the
same class $\overline{\alpha}\in \Z \pi_1(X^H(x),x)_{\phi'}$. We move
them into $X^{>H}(x)$ whenever possible. (We might need
to create a fixed point orbit in the inessential fixed point class
before~\cite[1.1]{wong91_location}.) We remove the inessential fixed
point orbits. We arrive at a map $h\simeq_G f$ such that
${N_G(f)}_{\overline{x}}=\# \{\text{fixed point orbits of }h_H(x)\}\geq
{M_G(f)}_{\overline{x}}$. Using Proposition~\ref{proplowerbound}, we
obtain equality.  

2)
Since $X$ is a cocompact smooth $G$-manifold, there is a map $f'$
which is $G$-homotopic to $f$ and which only has finitely
many fixed point orbits. We have a partial ordering on the
$\overline{y}\geq \overline{x}$ given by $\overline{y}\geq
\overline{z}\Leftrightarrow \Mor(z,y)\neq \emptyset$. We apply
Lemma~\ref{lemcoalesce} to $f'$ to coalesce fixed point orbits of the
same class, starting from the top. Note that when we remove fixed
point orbits, we can only move them up in this partial ordering. That
is why the definition has to be so complicated. We remove the inessential fixed point orbits. We are
left with one fixed point orbit for every essential class. 

We now look at a class $\mathcal{C}$ such that $N^G(f)_{\overline{x}}=\#
\mathcal{C}$, and we coalesce the essential fixed point orbits with
the corresponding classes appearing in $\mathcal{C}$. (If the
corresponding class in $\mathcal{C}$ is inessential, we might need
to create a fixed point orbit in this inessential fixed point class
before~\cite[1.1]{wong91_location}.) We obtain a map $h\simeq_G f$ which has exactly $\# \mathcal{C}$
fixed point orbits. Hence $N^G(f)_{\overline{x}}=\# \{\text{fixed
  point orbits of }h^H(x)\} \geq
M^G(f)_{\overline{x}}$. Using Proposition~\ref{proplowerbound}, we obtain equality.
\end{proof}

In general, it is not possible to find a map $h\simeq_G
f$ realizing all minima simultaneously. As an example, one can take $G=\Zz$ acting on $X=S^4$ as an
involution so that $X^{\Zz}=S^3$. One obtains
$M_G(\id_{S^4})_{\overline{x}}=0$ for all $\overline{x}\in \Is \Pi(\Zz,S^4)$, but the minimal number of fixed
points in the $G$-homotopy class of the identity $\id_{S^4}$ is equal
to $1$~\cite[Remark~3.4]{wong93}. In this example, the standard gap 
hypotheses are not satisfied. Other examples where the standard gap 
hypotheses do not hold and where the converse of the equivariant Lefschetz theorem is false are given in~\cite[Section~5]{ferrario99}.

\section{The $G$-Jiang Condition}\label{sec5}

In the non-equivariant case, we know that the generalized Lefschetz invariant is the right
element when looking for a precise count of fixed points. We read off
the Nielsen numbers from this invariant. In general, the Lefschetz
number contains too little information. But under certain conditions, we can
conclude facts about the Nielsen numbers from the Lefschetz numbers
directly, and thus obtain a converse of the Lefschetz fixed point
theorem.  

These conditions are called Jiang conditions. See Jiang~\cite[Definition~II.4.1]{jiang}, where one can find a thorough
treatment, and Brown~\cite[Chapter~VII]{brown}. The Jiang group is a subgroup
of $\pi_1(X,f(x))$~\cite[Definition~II.3.5]{jiang}. We generalize its definition to the equivariant case. 

\begin{Def}
Let $G$ be a discrete group, let $X$ be a finite proper $G$-CW-complex, and let
$f\colon X\to X$ be a $G$-equivariant endomorphism. Then a $G$-equivariant self-homotopy
$h\colon f\simeq_G f$ of $f$ determines a path $h(x,-)\in
\pi_1(X^H(x),f(x))$ for every $\overline{x}\in
\Is\Pi(G,X)$, with $x\colon G/H\to X$. Define the \emph{$G$-Jiang
  group} of $(X,f)$ to be 
\begin{align*}
J_G(X,f) & := \biggl\{ \sum_{\overline{x}\in\Pi(G,X)} \bigl[h(x,-)\bigr] \, \bigg| \, h\colon f
\simeq_G f \; G \text{-equivariant self-homotopy} \biggr\} \\
 & \leq \bigoplus_{\overline{x}\in\Pi(G,X)} \pi_1\bigl(X^H(x),f(x)\bigr),
\end{align*}
and define the \emph{$G$-Jiang group} of $X$ to be
\begin{align*}
J_G(X) & := \biggl\{ \sum_{\overline{x}\in\Pi(G,X)} \bigl[h(x,-)\bigr] \, \bigg| \, h\colon \id
    \simeq_G \id \; G \text{-equivariant self-homotopy} \biggr\} \\
 & \leq
\bigoplus_{\overline{x}\in\Pi(G,X)} \pi_1\bigl(X^H(x),x\bigr).
\end{align*}
\end{Def}

In the non-equivariant case, we know that the Jiang group $J(X,f,x)$
is a subgroup of the centralizer of $\pi_1(f,x)\bigl(\pi_1(X,x)\bigr)$
in $\pi_1\bigl(X,f(x)\bigr)$. In particular, $ J(X)\leq
Z\bigl(\pi_1(X,x)\bigr)$, where $Z\bigl(\pi_1(X,x)\bigr)$ denotes the center of
$\pi_1(X,x)$~\cite[Lemma~II.3.7]{jiang}. Furthermore, the
isomorphism $(f\circ w)_*\colon \pi_1\bigl(X,f(x_1)\bigr)\to \pi_1\bigl(X,f(x_0)\bigr)$ induced by a path~$w$ from
$x_0$ to $x_1$ 
induces an isomorphism $(f\circ w)_*\colon J(X,f,x_1)\to J(X,f,x_0)$ which does
not depend on the choice of $w$. So the definition does not depend on
the choice of the basepoint~\cite[Lemma~II.3.9]{jiang}. It is also known that $J(X)\leq
J(X,f)\leq \pi_1(X)$ for all $f$~\cite[Lemma~II.3.8]{jiang}. This leads to the consideration of spaces
with $J(X)=\pi_1(X)$ in the definition of a Jiang space. All these
lemmata also make sense in the
equivariant case. Thus we make the following definition.

\begin{Def}\label{defGJiangspace}
Let $G$ be a discrete group and let $X$ be a cocompact $G$-CW-complex. Then $X$
is called a \emph{$G$-Jiang space} if for all $\overline{x}\in
\Is\Pi(G,X)$ we have 
\[
{J_G(X)}_{\overline{x}}=\pi_1\bigl(X^H(x),x\bigr).
\]
\end{Def}

The group $J_G(X,f)$ acts on $\Lambda_G(X,f)$ as follows: If
$X^H(f(x))=X^H(x)$ and $X^H(f(x))=X^H(x)$, then
${J_G(X,f)}_{\overline{x}}$ acts on $\Z \pi_1(X^H(x),x)$. The element
$u=[h(x,-)] \in
{J_G(X,f)}_{\overline{x}}$ acts as composition with
$[w u w^{-1}] $, where $v=(\id,[w])\in \Mor(f(x),x)$. Since
$J_G(X,f)_{\overline{x}}$ is contained in the centralizer of
$\pi_1(f^H(x),x)\bigl(\pi_1(X^H(x),x)\bigr)$ in
$\pi_1\bigl(X^H(x),f(x)\bigr)$, this action induces an action on
$\Z\pi_1(X^H(x),x)_{\phi'}$ by composition, whence on $\Lambda_G(X,f)_{\overline{x}}$. Thus $J_G(X,f)$ acts
on $\Lambda_G(X,f)$, and by invariance of 
$\lambda_G(f)$ under homotopy
equivalence, we see that
\[
\lambda_G(f)\in \bigl(\Lambda_G(X,f)\bigr)^{J_G(X,f)}.
\]

Examples for $G$-Jiang spaces can be obtained from Jiang spaces. It is
known~\cite[Theorem~II.3.11]{jiang} that the class of Jiang spaces is
closed under homotopy equivalence and the topological product
operation and contains
\begin{itemize}
\item
simply connected spaces
\item
generalized lens spaces
\item
H-spaces
\item
homogeneous spaces of the form $A/A_0$ where $A$ is a topological
group and $A_0$ is a subgroup which is a connected compact Lie group.
\end{itemize}

Hence we obtain many examples of $G$-Jiang spaces using the following
proposition, analogous to Wong~\cite[Proposition~4.9]{wong93}.

\begin{Prop}
Let $G$ be a discrete group, and let $X$ be a free cocompact connected
proper
$G$-space. If $X/G$ is a Jiang space, then $X$ is a $G$-Jiang space. 
\end{Prop}

\begin{proof}
Since $X$ is connected and free, the set 
$\Is\Pi(G,X)$ consists of one element. Let $x$ be a basepoint of $X$. We
need to check that
${J_G(X)}_{\overline{x}}=\pi_1(X,x)$. Let $X\xrightarrow{p} X/G$ be the
projection. The Jiang subgroup of $X/G$ is given by 
\begin{eqnarray*}
J(X/G) & := & \Bigl\{\bigl[h(p(x),-)\bigr]\, \Big| \, h\colon \id_{X/G}\simeq \id_{X/G} \text{ self-homotopy} \Bigr\} \\
& \leq & \pi_1\bigl(X/G,p(x)\bigr).
\end{eqnarray*}
Let $\alpha\in \pi_1(X,x)$. Since $X\xrightarrow{p} X/G$ is a discrete cover,
$\widetilde{X}=\widetilde{X/G}$. There is a map
$p_\#\colon \pi_1(X,x)\to \pi_1(X/G,p(x))$ induced by the projection. Since $X/G$ is a Jiang space,
$J(X/G)=\pi_1(X/G,p(x))$, so there is a homotopy $h\colon \id_{X/G}\simeq \id_{X/G}$ such that $p_\#(\alpha)=
[h(p(x),-)]$. Because of the free and proper action of $G$ on $X$, this
homotopy $h$ can be lifted to a $G$-equivariant homotopy $h'\colon \id_X
\simeq_G \id_X$ such that $\alpha=[h'(x,-)]$. Thereby $\alpha\in
J_G(X)$. 
\end{proof}

\section[Converse of the Equivariant Lefschetz Fixed Point Theorem]{The Converse of the Equivariant Lefschetz Fixed Point Theorem}\label{sec6}

One can derive equivariant analogs of statements about Nielsen
numbers found in Jiang~\cite{jiang}, generalizing results of Wong~\cite{wong93} to infinite
discrete groups. In particular, if $X$ is a $G$-Jiang space, the
converse of the equivariant Lefschetz fixed point theorem holds. We
have the following theorem, compare Jiang~\cite[Theorem~II.4.1]{jiang}.

\begin{Thm}\label{thm4.10*}
Let $G$ be a discrete group, and let $X$ be a finite proper 
$G$-CW-complex which is a $G$-Jiang space. Then for any $G$-map
$f\colon X\to X$ and $\overline{x}\in \Is\Pi(G,X)$ with $x\colon
G/H\to X$ we have:
\begin{eqnarray*}
{L_G(f)}_{\overline{x}}=0 & \implies & \lambda_G(f)_{\overline{x}}=0 \text{ and }{N_G(f)}_{\overline{x}}=0,\\
{L_G(f)}_{\overline{x}}\neq 0 & \implies & \lambda_G(f)_{\overline{x}}\neq 0 \text{
  and }{N_G(f)}_{\overline{x}}=\# \bigl\{ \pi_1(X^H(x),x)_{\phi''}\bigr\}.
\end{eqnarray*}
Here $L_G(f)$ is the equivariant Lefschetz class~\cite[Definition~3.6]{lueck-rosenberg}, the equivariant analog of
the Lefschetz number.
\end{Thm}

\begin{proof}
Since $X$ is a $G$-Jiang space, the $G$-Jiang group $J_G(X)$ acts transitively on
$\pi_1(X^H(x),x)$ for all $\overline{x}\in \Is\Pi(G,X)$. This implies that 
\[
\lambda_G(f)_{\overline{x}} = \sum_{\overline{\alpha}}
  n_{\overline{\alpha}} \cdot \overline{\alpha} = n \cdot
  \sum_{\overline{\alpha}} \overline{\alpha}
\]
for some $n \in \Z$. This leads to  
$
{L_G(f)}_{\overline{x}} = n \cdot \# \bigl\{ \pi_1(X^H(x),x)_{\phi'} \bigr\}
$
by the augmentation map. We see that 
\begin{eqnarray*}
{L_G(f)}_{\overline{x}} =0 & \implies & n= 0 \\
& \implies & \lambda_G(f)_{\overline{x}} = 0 \\
& \implies & \nu_G(f)_{\overline{x}} = 0 \\
& \implies & {N_G(f)}_{\overline{x}}=0 , \\
{L_G(f)}_{\overline{x}}\neq 0 & \implies & n \neq 0 \\
& \implies & \lambda_G(f)_{\overline{x},\overline{\alpha}}\neq 0 \text{ for
  all }
\overline{\alpha}\in \Z(\pi_1(X^H(x),x))_{\phi'} \\
& \implies & \nu_G(f)_{\overline{x},\overline{\alpha}}\neq 0 \text{ for
  all }
\overline{\alpha}\in \Z(\pi_1(X^H(x),x))_{\phi''} \\
& \implies &
{N_G(f)}_{\overline{x}}=\#
\bigl\{ \pi_1(X^H(x),x)_{\phi''} \bigr\}. 
\end{eqnarray*}
\vskip-4.5ex
\end{proof}

The proof of Theorem~\ref{thm4.10*} already works if $J_G(X,f)$ acts transitively on every summand of
$\Lambda_G(X,f)$. We could have called $X$ a $G$-Jiang space if the
condition that $J_G(X,f)$ acts transitively on every summand of
$\Lambda_G(X,f)$ is satisfied. But this condition is less
tractable. It is implied by
${J_G(X,f)}_{\overline{x}}=\pi_1(X^H(x),f(x))$ for all $\overline{x}$,
which is implied by ${J_G(X)}_{\overline{x}}=\pi_1(X^H(x),x)$
for all $\overline{x}$.  

We now show that $f$ is $G$-homotopic to a fixed point free $G$-map if
the generalized equivariant Lefschetz invariant $\lambda_G(f)$ is zero. 

\begin{Thm}\label{thmlambdanulldannfixpunktfrei}
Let $G$ be a discrete group. Let $X$ be a cocompact proper smooth $G$-manifold
satisfying the standard gap hypotheses. Let $f\colon X\to X$ be a
$G$-equivariant endomorphism. Then the following holds:
\[
\text{If } \lambda_G(f)=0, \text{ then $f$ is $G$-homotopic to a fixed point
  free $G$-map.}
\]
\end{Thm}

\begin{proof}
If $\lambda_G(f)=0$, then $\ch_G(X,f)(\lambda_G(f))=0$, and therefore
we have ${N^G(f)}_{\overline{x}}=0$ for all $\overline{x}\in \Is\Pi(G,X)$. We
know from Theorem~\ref{thmdoubtful} that
$N^G(f)_{\overline{x}}=M^G(f)_{\overline{x}}=\min \{ \# \text{ fixed
  point orbits of }\varphi^H(x) \, | \, \varphi\simeq_G f\}$. In
particular, for $x\colon G/\{1\}\to X$ we obtain a map $\varphi$ such
that $\varphi^{\{1\}}(x)$ is fixed point free and $\varphi\simeq_G
f$. Thus we obtain our result on every connected component of $X$, and
combining these we arrive at a map $h\simeq_G f$ which is fixed point
free.
\end{proof}

These two theorems, Theorem~\ref{thm4.10*} and
Theorem~\ref{thmlambdanulldannfixpunktfrei}, combine to give the main
theorem of this paper, the converse of the equivariant Lefschetz fixed point theorem. 

\begin{Thm}\label{thmconvequL}
Let $G$ be a discrete group. Let $X$ be a cocompact proper smooth $G$-manifold
satisfying the standard gap hypotheses which is a
$G$-Jiang space. Let $f\colon X\to X$ be a
$G$-equivariant endomorphism. Then the following holds:
\[
\text{If } L_G(f)=0, \text{ then $f$ is $G$-homotopic to a fixed point
  free $G$-map.}
\]
\end{Thm}

\begin{proof}
We know that $L_G(f)=0$ means that ${L_G(f)}_{\overline{x}}=0$ for
all~$\overline{x}\in \Is\Pi(G,X)$. Since $X$ is a $G$-Jiang space, by
Theorem~\ref{thm4.10*} this implies that $\lambda_G(f)_{\overline{x}}=0$ for all
$\overline{x}\in \Is\Pi(G,X)$, so we have $\lambda_G(f)=0$. We
apply Theorem~\ref{thmlambdanulldannfixpunktfrei} to arrive at the
desired result.
\end{proof}

\begin{Rem}
As another corollary of Theorem~\ref{thmlambdanulldannfixpunktfrei}, we obtain:
If $G$ is a discrete group and $X$ is a cocompact proper smooth $G$-manifold
satisfying the standard gap hypotheses, then $\chi^G(X)=0$ implies
that the identity~$\id_X$ is $G$-homotopic to a fixed point free $G$-map. This was
already stated by L\"uck and
Rosenberg~\cite[Remark~6.8]{lueck-rosenberg}. Here $\chi^G(X)$ is the universal equivariant Euler
characteristic of $X$~\cite[Definition~6.1]{lueck-rosenberg} defined
by $\chi^G(X)_{\overline{x}}=\chi\bigl(WH_x\setminus X^H(x), WH_x\setminus
X^{>H}(x)\bigr)\in \Z$, we have $\chi^G(X)=L_G(\id_X)$. We calculate that 
\begin{eqnarray*}
\lambda_G(\id_X)_{\overline{x}} & = & \sum_{p\geq 0} (-1)^p \sum_{\Aut(x)
  \cdot e \in \atop \Aut(x)\setminus
  I_p(\widetilde{X^H(x)},\widetilde{X^{>H}(x)})}
\inc_{\phi}\bigl(\id_{\widetilde{X^H(x)}},e\bigr)\\ 
& = & \chi\bigl(WH_x\setminus X^H(x), WH_x\setminus X^{>H}(x)\bigr) \cdot
\overline{1}\in \Z \pi_1\bigl(X^H(x),x\bigr)_{\phi'}.
\end{eqnarray*} 
So we have $\chi^G(X)=0$ if and only if $\lambda_G(\id_X)=0$, and with
Theorem~\ref{thmlambdanulldannfixpunktfrei} we conclude that there is an
endomorphism $G$-homotopic to the identity which is fixed point free.
\end{Rem}

\bibliographystyle{alpha}

\def\cprime{$'$}

\end{document}